\documentclass[11pt,a4paper]{article}
\usepackage[cp1251]{inputenc}
\usepackage[english]{babel}
\usepackage{amsmath}
\usepackage{amssymb}
\usepackage{amsfonts}
\usepackage{authblk}

\title{"Painlev\'e 34" equation: equivalence test}
\author[1]{ Vera V. Kartak \thanks{kvera@mail.ru}}
\affil[1]{Laboratory "Group Analysis of Mathematical Models in Natural and Engineering Sciences", Ufa State Aviation Technical University,  K.Marx Str., 12, Ufa, 450000, Russia}

\affil[1]{Chair of Higher Algebra and Geometry, Faculty of Mathematics and Information Technologies, Bashkir State University, Z.Validi Str., 32, Ufa, 450076, Russia}

\begin{document}
\date{}
\maketitle 

{\small
\begin{quote}
\noindent{\bf Abstract. } We give the complete solution of the Equivalence Problem for ``Painlev\'e 34'' equation.
\medskip

{\bf Keywords:} {Invariant; Problem of Equivalence; Point transformation; Painlev\'e equation; Painlev\'e property}

{\bf 2000 Mathematics Subject Classification:} 53A55, 34A26, 34A34, 34C14, 34C20, 34C41 
\end{quote}
}

\section{Introduction. "Painlev\'e 34" equation}

At the beginning of XX century P. Painlev\'e and others studied
the following class of second order ODEs
$$
y''=F(x,y;y'),
$$
where the function $F$ is rational in $y'$ and analytic in $x$.
Their goal was to find all equations whose general solutions
 have no  movable critical singularities, i.e. have the Painlev\'e property.
 They solved this problem completely and found
 50 equations. Six of which were principally new -- {\it irreducible equations} -- (they did not allow reducing the order, and their solutions defined new special functions), they are currently called the {\it Painlev\'e equations} (PI-PVI equations),  see \cite{Gromak}, \cite{Kudryashov}.  
 In some books all forenamed 50 equations are named "Painlev\'e equation 1-50". The complete list of them is in books \cite{Gambier}, \cite{Ince}. 
 
A distinctive feature of  the "Painlev\'e 34"  equation is that its general solution and the PII solution 
\begin{equation}\label{p2}
PII:\qquad\tilde y''=2\tilde y^3+\tilde x\tilde y+a,\quad a=const
\end{equation}
are expressed one into the other explicitely using the B\"acklund transformation, see  \cite{Ince}, \cite{Fokas}, \cite{Yablonskii}, \cite{Kitaev}, \cite{Ohyama}. They can be written in the form of a Hamiltonian system of ordinary differential equations with one degree of freedom, see \cite{Suleimanov}.

Equation "Painlev\'e 34" from the book \cite{Ince} is 
\begin{equation}\label{ince34}
XXXIV.\qquad y''=\frac{{y'}^2}{2y}+4 a y^2-xy-\frac{1}{2y},\quad a=const\ne 0.
\end{equation}

In paper  \cite{Suleimanov}  this equation has some different form
\begin{equation}\label{sul34}
y''=\frac{y'^2}{2y}-2y^2-xy-\frac{(\alpha\pm 1/2)^2}{2y},\qquad \alpha= const.
\end{equation}

Note, that equation "Painlev\'e 34" plays important role in the description of multi-ion electro-diffusion models, see \cite{Rogers}.

In paper \cite{Hietarinta} was first stated the problem of deriving syzygies (relationships between the invariants) for every equation from the list of Gambier \cite{Ince} i.e. for every "Painlev\'e equation 1-50". This work was continued in the recent paper \cite{Bagderina1}, where was 
found some syzygies for the equation from the list of Gambier including equation (\ref{ince34}).

The aim of this paper is constructing the equivalence test -- necessary and sufficient conditions written in terms of invariants checking the equivalence of some equation (\ref{eq}) to the "Painlev\'e 34"
equation (\ref{sul34}) under the general point transformations
\begin{equation}\label{zam}
\tilde x=\tilde x(x,y),\quad \tilde y=\tilde y(x,y).
\end{equation}

\section{Implementation of the classification and calculation of the invariants}

Equation (\ref{sul34}) is from the
 following class of the second order ODE's
\begin{equation}\label{eq}
y''=P(x,y)+3\,Q(x,y)y'+3\,R(x,y)y^{\prime 2}+S(x,y)y^{\prime 3},
\end{equation}
that is the closed  under the general point transformations (\ref{zam}).

In the set of papers \cite{Sharipov1, Sharipov2, Sharipov3}, review see \cite{Kartak},
Ruslan Sharipov 
succeeded to construct the system of  (pseudo)invariants which he 
calculated explicitly in the terms of the coefficients of equations (\ref{eq}).
On the basis of this system he classified equations (\ref{eq}).
In the present paper we use this classification for solving the equivalence problem of equation (\ref{sul34}).

{\bf Step 1.} At first we write equation (\ref{sul34}) in the form
\begin{equation}\label{sul34b}
y''=\frac{y'^2}{2y}-2y^2-xy-\frac{{\beta}^2}{2y},\quad \beta=\alpha\pm 1/2.
\end{equation}
Equation (\ref{sul34b}) has the form (\ref{eq}) with the  coefficients
$$
P(x,y)=-2y^2-xy-\frac{{\beta}^2}{2y},\quad Q(x,y)=0,\quad R(x,y)=\frac{1}{6y},\quad S(x,y)=0.
$$

{\bf Step 2.} Then calculate the basic objects characterizing the equation (\ref{sul34b}). 
Details are in the papers \cite{Sharipov1, Sharipov2, Sharipov3},  \cite{Kartak}.

\textit{Pseudotensorial field of weight $ m $ and  valence $ (r, s) $} is an indexed set  transformed  under change of variables (\ref{zam}) by the rule
$$
F^{i_1\dots i_r}_{j_1\dots j_s}=(\det T)^m{\sum_{p_1\dots p_r}}
{\sum_{q_1\dots q_s}} S^{i_1}_{p_1}\dots S^{i_r}_{p_r}T^{q_1}_{j_1}
\dots T^{q_s}_{j_s} \tilde F^{p_1\dots p_r}_{q_1\dots q_s},
$$
where $S$ and $T$ are direct and inverse  transformations matrices for  (\ref{zam}).

The first pseudovectorial field ${\boldsymbol \alpha}$ associated with equation (\ref{eq}) 
has weight 2 and the components $\alpha^1=B$, $\alpha^2=-A$, where
\begin{equation}\label{alpha}
\aligned 
A&=P_{ 0.2}-2Q_{ 1.1}+R_{ 2.0}+ 2PS_{ 1.0}+\\
&+SP_{1.0}-3PR_{ 0.1}-3RP_{ 0.1} -3QR_{ 1.0} +6QQ_{ 0.1}, \\
 B&=S_{ 2.0}-2R_{ 1.1}+Q_{ 0.2}- 2SP_{0.1}-\\
&-PS_{ 0.1}+3SQ_{ 1.0}+3QS_{ 1.0}+ 3RQ_{ 0.1}-6RR_{ 1.0}.
\endaligned
\end{equation}
We can check this fact applying the direct symbolic calculations.
Hereinafter  symbol $K_{i.j}$ denotes the partial differentiation:   $K_{i.j}={\partial ^{i+j}K}/{\partial x^i\partial y^j}.$

The second pseudovectorial field $\boldsymbol{\beta}$ has  weight 4 and the components $\beta^1=G$, $\beta^2=H$, where
$$
\aligned
 G&=-BB_{ 1.0}-3AB_{ 0.1}+4BA_{ 0.1}+
3SA^2-6RBA+3QB^2,\\
H&=-AA_{ 0.1}-3BA_{ 1.0}+4AB_{ 1.0}-
3PB^2+6QAB-3RA^2.
\endaligned
$$
Their scalar product (using the skew-symmetric Gramian matrix) denoting the 
pseudoinvariant $F$  by the formula
\begin{equation}\label{F}
3F^5=AG+BH.
\end{equation}

For the equation (\ref{sul34b}) ${\boldsymbol \alpha}$  from (\ref{alpha}) and $F$ from ({\ref{F}}) are equal to
$$
A=-3-\frac{3{\beta}^2}{8y^3},\quad B=0,\quad F=0.
$$

{\bf Step 3.} As if for the equation (\ref{sul34b}) the conditions $F=0$, but $A\ne 0$ or $B\ne 0$
are true, it relates to the {\it Case of intermediate degeneration}, for the details see \cite{Sharipov3}. 
In this case we can calculate another important pseudoinvariants $\Omega$ of weight 1, $N$ of weight 2 and $M$ of weight 4 using explicit formulas, that are different in the cases $A\ne 0$ or $B\ne 0$.

As $A\ne 0$, the explicit formula for pseudoinvariant $\Omega$  reads as
\begin{equation}\label{Omega1}
\aligned \Omega &=\frac {2BA_{ 1.0}(BP+ A_{ 1.0})}{A^3}- \frac
{(2B_{ 1.0}+3BQ)A_{ 1.0}}{A^2}+\frac {(A_{ 0.1}-2B_{1.0})BP}{A^2}-\\
&- \frac {BA_{ 2.0}+B^2 P_{ 1.0}}{A^2}+ \frac
{B_{ 2.0}}A+\frac {3B_{ 1.0}Q+3BQ_{ 1.0}- B_{ 0.1}P-BP_{
0.1}}{A}+\\
&+Q_{ 0.1}- 2R_{ 1.0}.
\endaligned
\end{equation}
And in the case $B\ne 0$ the similar formula is
\begin{equation}\label{Omega2}
\aligned \Omega &=\frac {2AB_{ 0.1}(AS- B_{ 0.1})}{B^3}- \frac
{(2A_{ 0.1}-3AR)B_{ 0.1}}{B^2}+\frac {(B_{ 1.0}-2A_{
0.1})AS}{B^2}+ \\
&+\frac {AB_{ 0.2}-A^2 S_{ 0.1}}{B^2}- \frac
{A_{ 0.2}}B+\frac {3A_{ 0.1}R+3AR_{ 0.1}- A_{ 1.0}S-AS_{
1.0}}{B}+\\
&+R_{ 1.0}- 2Q_{ 0.1}.
\endaligned
\end{equation}
In the cases $A\ne 0$ and $B\ne 0$ the pseudoinvariant $N$  is given by the formulas 
\begin{equation}\label{N}
 N =-\frac H{3A}, \qquad\qquad N =\frac G{3B}.
\end{equation}
The pseudoinvariant $M$  in the case $A\ne 0$ reads as
\begin{equation}\label{M1}
\aligned
M=&-\frac {12BN(BP+A_{ 1.0})}{5A}+\frac {24}5BNQ+\frac
65NB_{ 1.0}+\frac 65NA_{ 0.1}-\\
&-AN_{ 0.1}+BN_{ 1.0}- \frac {12}5ANR.
\endaligned
\end{equation}
And in the case $B\ne 0$ reads as
\begin{equation}\label{M2}
\aligned
M=&-\frac {12AN(AS-B_{ 0.1})}{5B}+\frac {24}5ANR-
 \frac
65NA_{ 0.1}-\frac 65NB_{ 1.0}+\\
&+BN_{ 1.0}-AN_{ 0.1}- \frac {12}5 BNQ.
\endaligned
\end{equation}

For the equation (\ref{sul34b}) pseudoinvariants $\Omega$ from (\ref{Omega1}), $N$ from (\ref{N})  and $M$ from (\ref{M1}) are equal to
$$
\quad\Omega=0,\qquad N=\frac{5{\beta}^2}{4y^4}-\frac{1}{2y},\qquad M=\frac{9}{10y^2}-\frac{63{\beta}^2}{4y^5}.
$$

{\bf Step 4.} 
It is easy to see that the pseudoinvariant $M$ given by (\ref{M1}), (\ref{M2}) for the equation (\ref{sul34b}) is not vanishing.

As if $M\ne 0$ for the equation (\ref{sul34b}), then it relates to the {\it First case of intermediated degeneration}, see \cite{Sharipov3}. In this case the basic invariants are 
\begin{equation}\label{inv1}
I_1=\frac{M}{N^2},\qquad I_2=\frac{\Omega^2}{N},\qquad I_3=\frac{\hat\Gamma^1_{22}}{M},\quad\text{where}
\end{equation}
$$
\aligned
\hat\Gamma^1_{22}=&\frac {\gamma^1\gamma^2(\gamma^1_{
1.0}- \gamma^2_{ 0.1})}{M}+ \frac {(\gamma^2)^2\gamma^1_{ 0.1}-
(\gamma^1)^2\gamma^2_{ 1.0}}M+\\
&+\frac
{P(\gamma^1)^3+3Q(\gamma^1)^2\gamma^2+3R\gamma^1(\gamma^2)^2+
S(\gamma^2)^3}M.
\endaligned
$$
Here ${\boldsymbol \gamma}$ is a new pseudovectorial field of weight 3 associated with equation (\ref{eq}) relating to the First case of intermediate degeneration

As $A\ne 0$, the components of the pseudovectorial field ${\boldsymbol \gamma}$ reads as
\begin{equation}\label{gamma1}
\aligned
 \gamma^1=&-\frac {6BN(BP+A_{ 1.0})}{5A^2}+
\frac {18NBQ}{5A}+
\frac {6N(B_{ 1.0}+A_{ 0.1})}{5A} -\\
&-N_{ 0.1}-\frac
{12}5NR-2\Omega B,\\
\gamma^2=&-\frac {6N(BP+A_{ 1.0})}{5A}+N_{ 1.0}+\frac 65NQ+ 2\Omega
A.
\endaligned
\end{equation}

As $B\ne 0$,  reads as
\begin{equation}\label{gamma2}
\aligned
\gamma^1=&-\frac {6N(AN-B_{ 0.1})}{5B}-N_{ 0.1} +\frac 65NR-2\Omega B,\\
 \gamma^2=&-\frac {6AN(AS-B_{ 0.1})}{5B^2}+
\frac {18NAR}{5B}-
\frac {6N(A_{ 0.1}+B_{ 1.0})}{5B} +\\
&+N_{ 1.0}-\frac {12}5NQ+2\Omega A.
\endaligned
\end{equation}

The additional invariants are computed by differentiating  the basic invariants (\ref{inv1}) along  pseudo\-vectorial fields
${\boldsymbol \alpha}$ from (\ref{alpha}) and ${\boldsymbol \gamma}$ from (\ref{gamma1}), (\ref{gamma2})

\begin{equation}\label{inv1-1}
\begin{aligned}
I_4&=\frac{B(I_1)'_x-A(I_1)'_y}{N},\qquad\qquad I_6=\frac{B(I_3)'_x-A(I_3)'_y}{N},\\
I_7&=\frac{(\gamma^1(I_1)'_x+\gamma^2(I_1)'_y)^2}{N^3},\qquad I_9=\frac{(\gamma^1(I_3)'_x+\gamma^2(I_3)'_y)^2}{N^3},\\
 I_{15}&=\frac{(\gamma^1(I_6)'_x+\gamma^2(I_6)'_y)^2}{N^3},\qquad I_{21}=\frac{(\gamma^1(I_9)'_x+\gamma^2(I_9)'_y)^2}{N^3}.
\end{aligned}
\end{equation}

For the equation (\ref{sul34b}),  invariants $I_1$, $I_2$ from (\ref{inv1}) and $I_7$ from (\ref{inv1-1}) are  
\begin{equation}\label{inv1p34}
I_1= -\frac{36}{5}\frac{y^3(35{\beta}^2-2y^3)}{(5{\beta}^2-2y^3)^2},\qquad I_2=0,\qquad I_7=0.
\end{equation}
As we can see, the invariant $I_7$ is vanishing, so  the equation (\ref{sul34b}) relating to the {\it  Case 1.4 of intermediate degeneration}, for details see \cite{Kartak}.

\section{Equivalence test}

It may be two different possibilities, $I_1$ given by (\ref{inv1p34}) is a constant or not.

\subsection{Case $I_1=const$}
Equations (\ref{eq}) relating to the {\it  First case  of intermediate degeneration} with
the conditions $I_1=const$, $I_2=0$ from (\ref{inv1}) were  described in paper \cite{Kartak2}.

Let us represent formula for $I_1$ from (\ref{inv1p34}) in the form
$$
I_1=\frac{18}{5}-\frac{90{\beta}^2(2y^3+{\beta}^2)}{(5{\beta}^2-2y^3)^2}.
$$
It is not difficult to see that the only way $I_1$ to be a constant is $\beta=0$.

Then we calculate the invariant $I_3$ from (\ref{inv1}), the additional invariants $I_6$, $I_9$ from (\ref{inv1-1}) and a new invariant $J$, where
\begin{equation}\label{j}
J=\frac {4+10I_6-60I_3}{50\sqrt{I_9}}.
\end{equation}

For the equation (\ref{sul34b}) with the zero parameter $\beta$  these invariants are 
$$
I_3=\frac1{30}\frac{2y+x}{y},\quad I_6=\frac 15\frac xy,\quad I_9=-\frac{1}{1250}\frac 1{y^3},\quad I_{21}=0,\quad J=0.
$$

In papers \cite{Kartak}- \cite{Kartak2} the following Theorem was proved. (In paper \cite{Kartak} the condition $I_{21}=0$ is unfortunately missed.)

{\bf Theorem 1.} {\it Equation (\ref{eq}) is equivalent to Painleve II equation (\ref{p2}) with the parameter  $a=\pm J$  if and only if the following conditions hold: 
\begin{enumerate}
\item  equation corresponds to the Case of intermediate degeneration: $A\ne 0$ or $B\ne 0$ in (\ref{alpha}), but $F=0$ in (\ref{F}); 
\item equation corresponds to the First case  of intermediate degeneration: $M\ne 0$ in (\ref{M1}), (\ref{M2}), $\Omega=0$ in (\ref{Omega1}), (\ref{Omega2});
\item  $I_1=18/5$ in (\ref{inv1}), $I_9\ne 0$, $I_{21}=0$ in (\ref{inv1-1}), invariant $J=const$ in (\ref{j}). Among the invariants $I_3$, $I_6$ and $I_9$ from (\ref{inv1}), (\ref{inv1-1}) one can find two functionally independent.
\end{enumerate}
 The  invariant point transformation is 
$$
\tilde y= \frac{1}{\sqrt[6]{2500 I_9}},\qquad \tilde x=\frac{5I_6}{\sqrt[6]{2500I_9}}-\frac 32 J\sqrt[6]{2500 I_9}. $$
}

For the equation (\ref{sul34b}) with zero parameter $\beta$ all conditions of Theorem 1 are hold.
So it is equivalent to Painlev\'e II equation (\ref{p2}) with zero parameter $a$. The corresponding 
change of variables 
$x=-\sqrt[3]{2}\tilde x,$ $ y=-\sqrt[3]{2}{\tilde y}^2$
transforms equation (\ref{sul34b}) (that is written in variables $(x,y)$) into  equation (\ref{p2}) (that is written into variables $(\tilde x,\tilde y)$) with the parameter $a=0.$

\subsection{Case $I_1\ne const$}

It was proved above that in this case the parameter $\beta$ is not vanishing.
Let us make the following change of variables,
$$
y = {\tilde y}^{1/3} {\beta}^{2/3},\qquad x = {\tilde x} {\beta}^{2/3},
$$
then the equation (\ref{sul34b}) takes the form
\begin{equation}\label{sul34bb}
{\tilde y''}=\frac {5{\tilde y}^{\prime 2}}{6\tilde y}-{\beta}^2{\tilde y}^{1/3}\left(6{\tilde y}+3{\tilde x}{\tilde y}^{2/3}+\frac 32\right).
\end{equation}

To simplify the notation below we do not write the tildes over the variables $x$ and $y$ in the equation (\ref{sul34bb}). This equation also has form (\ref{eq}) with the coefficients
$$
 P=-{\beta}^2{ y}^{1/3}\left(6{ y}+3{ x}{ y}^{2/3}+\frac 32\right),\quad Q=0,\quad
   R=\frac 5{18 y},\quad S=0.
$$

Let us calculate invariants $I_1$, $I_3$ from (\ref{inv1}),   $I_4$, $I_7$, $I_9$, $I_{15}$, $I_{21}$ from (\ref{inv1-1})
for the equation (\ref{sul34bb})
\begin{equation}\label{inv_p34}
\aligned
I_1=& \frac{36}5 \frac{y(2y-35)}{(2y-5)^2},\qquad  I_3= \frac{y(4y+2xy^{2/3}+1)(2y-35)}{15(2y+1)^3},\\
I_4 =& -\frac{3240(2y+7)y(2y+1)}{(2y-5)^4},\qquad I_9=-\frac{64}{625}\frac{y^6(2y-35)^4}{{\beta}^2(2y+1)^8(2y-5)^3},\\
I_6=&\frac{4y^{5/3}x(4y^2-296y+175)-6y(300y^2-136y-35)}{5(2y-5)(2y+1)^3},\quad 
 I_7=0,\\
I_{15}=&-\frac{2304y^6(4y^2-296y+175)^2(2y-35)^2}{625(2y-5)^2(2y+1)^8{\beta}^2},\qquad I_{21}=0.
\endaligned
\end{equation}
Then we regard the symbols $I_1$ and $I_4$ as the parameters in order to convert formulas for
$I_1$ and $I_4$ from (\ref{inv_p34}) into polynomials. We get two polynomials depending only on the variable $y$
$$
\mathbf{P_1}=36y(-35+2y)-5I_1(2y-5)^2,\qquad \mathbf{P_2}=3240(2y+7)y(2y+1)+I_4(2y-5)^4.
$$

By implementation of Buchberger's  algorithm, see \cite{Cox}, we reduce polynomials $\mathbf{P_1}$ and $\mathbf{P_2}$ with respect to the variable $y$. 
We obtain a new invariant $K$, that is vanishing for the equation (\ref{sul34bb})
\begin{equation}\label{k}
\aligned
K = &500 I_1^4-7275 I_1^3+500 I_4 I_1^2+32940 I_1^2-\\
&-5475 I_4 I_1-47628 I_1+125 I_4^2+13230 I_4=0.
\endaligned
\end{equation}

From the next-to-last step of Buchberger's algorithm, we get a formula for the variable $y$ in terms of invariants
$$
y=\frac{125(2322I_1+3 I_4+20 I_4I_1-915I_1^2+75I_1^3)}{2(-1469664+1250I_4I_1-13875I_4+691470I_1-90825I_1^2+3375I_1^3)}.
$$
The variable $x$ we find using the formula of $I_3$ from (\ref{inv_p34})
$$
x=\frac {(120 I_3-8)y^3+(138+180I_3)y^2+(90I_3+35)y+15I_3}{2y^{5/3}(2y-35)}.
$$

The parameter $\beta^2$ we find using the formula of $I_9$ from (\ref{inv_p34})
$$
\beta^2=-\frac{64}{625}\frac{y^6(2y-35)^4}{I_9(2y+1)^8(2y-5)^3}.
$$

The invariants $I_6$, $I_3$, $I_1$ and $I_{15}$, $I_9$  are related by the formulas
$$
\aligned
 I_6=&\frac{6(4 y^2-296 y+175)I_3}{(2 y-5)(2 y-35)}-
\frac{(2 y-1)(2 y-5)I_1}{9(2 y+1)^2},\\
I_{15}=&\frac{36(4 y^2-296 y+175)^2I_9}{(2 y-35)^2(2 y -5)^2}.
\endaligned
$$

{\bf Theorem 2. }{\it Equation (\ref{eq}) is equivalent to the "Painlev\'e 34" equation (\ref{sul34bb})  with the parameter $\beta\ne 0$ if and only if the following conditions hold:
\begin{enumerate}
\item the equation corresponds to the Case of intermediate degeneration: $A\ne 0$ or $B\ne 0$ from (\ref{alpha}), but $F=0$ from (\ref{F});
\item the equation corresponds to the Case 1.4 of intermediate degeneration: $M\ne 0$ from (\ref{M1}), (\ref{M2}), $I_2=0$ from (\ref{inv1}), $I_7=0$, $I_{21}=0$  from (\ref{inv1-1});  
\item the invariant $K=0$ from (\ref{k});
\item there exists a non-degenerate invariant change of variables
that connects equations (\ref{eq}) and (\ref{sul34bb}) 
\begin{equation}\label{y}
\aligned
&\tilde y=\frac{125}{2}\cdot \\
&\cdot\frac{(3+20I_1)I_4+3I_1(5I_1-18)(5I_1-43)}{125(10 I_1-111)I_4+3(5I_1-18)(225I_1^2-5245I_1+27216)}
\endaligned
\end{equation}
\begin{equation}\label{x}
\tilde x=\frac {(120 I_3-8){\tilde y}^3+(138+180I_3){\tilde y}^2+(90I_3+35){\tilde y}+15I_3}{2{\tilde y}^{5/3}(2{\tilde y}-35)},
\end{equation}
\item  the following invariant is a constant
\begin{equation}\label{b}
\beta^2=-\frac{64}{625}\frac{{\tilde y}^6(2\tilde y-35)^4}{I_9(2{\tilde y}+1)^8(2{\tilde y}-5)^3}.
\end{equation}
\item invariants $K_1=0$ and $K_2=0$, where
\begin{equation}\label{k1}
\aligned
K_1 = & I_6-\frac{6(4\tilde y^2-296\tilde y+175)I_3}{(2\tilde y-5)(2\tilde y-35)}+
\frac{(2\tilde y-1)(2\tilde y-5)I_1}{9(2\tilde y+1)^2},\\
K_2=&\frac{I_{15}}{I_9}-\frac{36(4\tilde y^2-296\tilde y+175)^2}{(2\tilde y-35)^2(2\tilde y -5)^2}.
\endaligned
\end{equation}
Here in the formulas (\ref{x}), (\ref{b}), (\ref{k1}) we should substitute the expression of $\tilde y$ via the invariants $I_1$ and $I_4$ from (\ref{y}).
\end{enumerate}
}

{\bf Example 1. }{  Let us return to equation (\ref{ince34}). All
conditions of Theorem 2 are true. The point transformation
$$
\tilde y=-2ay^3,\qquad \tilde x=\frac x{(2a)^{2/3}},\qquad a\ne 0
$$
transforms equation (\ref{ince34}) (that is written in variables $(x,y)$) into  equation (\ref{sul34bb}) (that is written into variables $(\tilde x,\tilde y)$) with the parameter $\beta^2=4a^2.$

Let's note, that in the case $a=0$ equation (\ref{ince34}) is equivalent to 
$y''=y^{-3},$ see \cite{Kartak} for the details.
}

{\bf Example 2. } { Equation (\ref{sul34}) is not equivalent to Painlev\'e IV equation
$$
y''=\frac{{y'}^2}{2y}+\frac{3y^3}{2}+4xy^2+2(x^2-\alpha)y-\frac{{\beta}^3}{2y},\quad \alpha,\, \beta=const.
$$
Indeed, for the equation PIV the invariants $I_7\ne 0$. See \cite{Kartak3}.
}

{\bf Example 3. } { Equation describing 3-ion case (3a) from \cite{Rogers}
$$
w''-\frac{w'^2}{2w}+\nu_1^2\left(-2k_1w^2-(Cx+K)w+\frac{k_2}{w}\right)=0,\quad \;\nu_1,\, k_1,\, k_2,\, C,\, K=const
$$
is equivalent to "Painlev\'e 34" equation (\ref{sul34bb}) with the parameter $\beta^2=2k_2\nu_1^2k_1^2/C^2$ if $\nu_1\ne 0,\; k_1\ne 0,\; k_2\ne 0,\; C\ne 0$. All conditions of Theorem 2 are true. The following point transformation
$$
\tilde y=-\frac{k_1w^3}{2k_2},\quad \tilde x=-\frac{Cx+K}{2^{1/3}k_1^{2/3}k_2^{1/3}}
$$
transforms this equation (that is written in variables $(x,w)$) into the equation (\ref{sul34bb}) (that is written into variables $(\tilde x,\tilde y)$).

And it is equivalent to Painlev\'e II equation (\ref{p2}) with the parameter $a=0$ if $\nu_1\ne 0,\; k_1\ne 0,\; k_2= 0,\; C\ne 0$. All conditions of Theorem 1 are true. The following point transformation
$$
\tilde y=\frac{\sqrt[3]{\nu_1}\sqrt{k_1w}}{\sqrt[6] 2\sqrt[3]{C}},\quad \tilde x=\frac{(Cx+K)\sqrt[3]{\nu_1}}{\sqrt[6] 2\sqrt[3]{C}\sqrt{k_1w}}
$$
transforms this equation (that is written in variables $(x,w)$) into the equation (\ref{p2}) (that is written into variables $(\tilde x,\tilde y)$).

}

{\bf Example 4. }{Equation describing 3-ion case (3b) from \cite{Rogers}
$$
\aligned
\left(w+\frac{Cx+K}{k_1}\right)w''&-\frac{w'^2}{2}-\frac{Cw'}{k_1}-2k_1\nu_1^2 w^3-\\
&-4\nu_1^2(Cx+K)w^2-2\nu_1^2(Cx+K)^2\frac{w}{k_1}=0,
\endaligned
$$
$\nu_1,\, k_1,\,\, C,\, K=const $ is equivalent to "Painlev\'e 34" equation (\ref{sul34bb}) with $\beta^2=1/4$ if $\nu_1\ne 0,\; k_1\ne 0,\;  C\ne 0$. All conditions of Theorem 2 are true. The following point transformation
$$
\tilde y=-\frac{\nu_1^2(k_1w+Cx+K)^3}{C^2},\quad \tilde x=\frac{2{\nu_1}^{2/3}(Cx+K)}{C^{2/3}}
$$
transforms this equation (that is written in variables $(x,w)$) into  the equation (\ref{sul34bb}) (that is written into variables $(\tilde x,\tilde y)$).
}

{\bf Example 5. } Reduction of Nonlinear Schro\"edinger equation, see \cite{Bagderina}. Here function $V(x)$ is a potential.
$$
y''=V(x)y-y^3+\frac{k^2}{y^3},\quad k=const.
$$
If $k\ne 0$, invariants are
$$
\aligned
A&=-6y+\frac{12k^2}{y^5},\quad B=0,\quad F=0,\quad M=\frac{72}{5}+\frac{1008k^2}{y^6},\\
I_2&=0,\quad I_7=0,\quad I_{21}=-\frac{288V'^2(x)y^{88}(y^6+70k^2)^{10}V''^2(x)}{9765625(y^6-2k^2)^{18}(y^6+10k^2)^9}=0,\\
K&=0,\quad  \beta^2=-\frac{k^2}{V'^2(x)}=const,\quad K_1=0,\quad K_2=0.
\endaligned
$$
So, all conditions of Theorem 2 are true if $k\ne 0$ and $V(x)$ is a certain linear function, then this equation is equivalent to "Painlev\'e 34" equation (\ref{sul34bb}). The following point transformation
$$
\tilde y=-\frac{y^6}{4k^2},\quad \tilde x=\sqrt[3]{2} k V(x)
$$
transforms this equation (that is written in variables $(x,y)$) into  the equation (\ref{sul34bb}) (that is written into variables $(\tilde x,\tilde y)$).

If $k=0$ then the invariants are
$$
\aligned
I_1&=\frac{18}{5},\quad I_3=\frac 1{15}-\frac{V(x)}{15y^2},\quad I_6=-\frac{2V(x)}{5y^2},\quad I_9=-\frac{2V'^2(x)}{625y^6},\\
I_{21}&=-\frac{288V'^2(x)V''^2(x)}{9765625y^{14}}=0,\quad J=0.
\endaligned$$
All conditions of Theorem 1 are true if $V(x)$ is a certain linear function,
$V(x)\ne const.$
The following linear point transformation
$$
\tilde y=\frac{\sqrt[6]{-1}y}{\sqrt{2}\sqrt[3]{V'(x)}},\qquad
\tilde x=\frac{V(x)}{\sqrt[3]{V'^2(x)}}
$$
transforms this equation (that is written in variables $(x,y)$) into  the equation Painlev\'e II (\ref{p2}) with zero parameter $a$ (that is written into variables $(\tilde x,\tilde y)$).

\subsection*{Acknowledgements}

I am grateful to Prof.\ Bulat Ir. Suleimanov for his useful remarks.

The work is partially supported by the Government of Russian Federation through Resolution No. 220, Agreement No. 11.G34.31.0042.


\begin{thebibliography}{99}


\bibitem{Gromak}	V. I. Gromak, N. A. Lukashevich,  Analytical Properties of Solution of Painlev\'e Equations, Izdatel'stvo Univer., Minsk, 1990, [in Russian].

\bibitem{Kudryashov} N. A. Kudryashov, Analytic Theory of Nonlinear Differential Equations,  IKI, Moscow, 2004, [in Russian].

\bibitem{Gambier} B. Gambier, Sur les equations differentielles du second ordre et du premier degre dont l'integrale generale est a points critiques fixes, \textit{Acta Math.} \textbf{33} (1910), 1-55.

\bibitem{Ince} E.L. Ince, Ordinary differential equations, Courier Dover Publications, 1956.

\bibitem{Fokas} A.S. Fokas, M.J. Ablowitz, On a Unified Approach to Transformations and Elementary Solutions of Painleve Equations
{http://dx.doi.org/10.1063/1.525260},
{\textit{J. Math. Ph.}} \textbf{23}(11) (1982) 2033-2042.

\bibitem{Yablonskii} A.I. Yablonskii, On rational solutions of the second Painleve equation, {\textit{ Vesti Akad. Navuk. BSSR Ser. Fiz. Tkh. Nauk.}} \textbf{3} (1959)
30–35 [in Russian].

\bibitem{Kitaev} L.A. Bordag, A.B. Kitaev, On the Connection of the Second Kind Painlev\'e Equation with Nonlinear Evolution Equations, {\textit{Communications JINR}}, P5-87-208, Dubna, 1987, [in Russian].

\bibitem{Ohyama} Y. Ohyama, S. Okumura, A coalescent diagram of the Painlev\'e equations from the viewpoint of isomonodromic deformations, {http://dx.doi.org/10.1088/0305-4470/39/39/S08}m
{\textit{ J. Phys. A: Math. Theor.}} \textbf{39}(39) (2006), 12129-12151.
 
\bibitem{Suleimanov} B. I. Suleimanov, “Quantizations” of the second Painlev\'e equation and the problem of the equivalence of its L-A pairs, 
{http://dx.doi.org/10.1007/s11232-008-0106-8},
{\textit{Theoretical and Mathematical Physics}}  \textbf{156}(3) (2008), 1280-1291.


\bibitem{Rogers} R. Conte, C. Rogers and W. K. Schief, Painlev\'e structure of a multi-ion electrodiffusion system, 
{http://dx.doi.org/10.1088/1751-8113/40/48/F01}
{\textit{ J. Phys. A: Math. Theor.}} \textbf{40} (2007), F1031-F1040.

\bibitem{Hietarinta} J. Hietarinta, V. Dryuma, Is my ODE a Painlev\'e equation in disguise?
{http://dx.doi.org/10.2991/jnmp.2002.9.s1.6},
	 {\textit{Journal of Nonlin. Math. Phys.}} \textbf{9}(1) (2002), 67-74.

\bibitem{Sharipov1} V.~V. Dmitrieva, R.~A. Sharipov,  On the
point transformations for the second order differential equations, {\it Electronic archive at LANL}, (1997) {solv-int/9703003}. 


\bibitem{Sharipov2} R.~A. Sharipov,  On the point transformations
for the equation $y''=P+3\,Q\,y'+3\,R\,{y'}^2+S\,{y'}^3$,
{\it Electronic archive at LANL}, (1997), {solv-int/9706003}.


\bibitem{Sharipov3} R.~A. Sharipov,  Effective procedure of point
classification for the equations $y''=P+3\,Q\,y'+3\,R\,{y'}^2
+S\,{y'}^3$, {\it Electronic archive at LANL}, (1998), {Math\.DG/9802027}.

\bibitem{Kartak} V. V. Kartak, 
Point classification of the second order ODEs and its application to Painlev\'e equations
\textit{  Journal of Nonlinear Math. Physics} 
{\bf 20}(supp. 1) (2013), 110-129. DOI:10.1080/14029251.2013.862438;
{\it Electronic archive at LANL} (2012), {1204.0174}.

\bibitem{Kartak1} V.~V.~Kartak,  Explicit solution of the equivalence problem for certain Painleve equations, {\textit{ Ufimskii Math. Journal}} 
\textbf {1}(3) (2009), 46--56, [in Russian].

\bibitem{Kartak2} V.~V.~Kartak,  
Equivalence classes of the second order ODEs
with the constant Cartan invariant, \textit{ Journal of Nonlinear Math. Physics} 
\textbf {18}(4) (2011), 613--640. 

\bibitem{Cox} David Cox, John Little, and Donal O'Shea,  { Ideals, Varieties, and Algorithms: An Introduction to Computational Algebraic Geometry and Commutative Algebra}, Springer, 1997.

\bibitem{Kartak3} V.~V.~Kartak, Solution of the equivalence problem for the Painleve IV equation, {http://dx.doi.org/10.1007/s11232-012-0132-4}, {\textit{  Theoretical and Mathematical Physics}}, \textbf{173}(2) (2012), 1541-1564. 

\bibitem{Bagderina} Yu. Yu. Bagderina, Equivalence of ordinary differential equations $y'' = R(x, y)y'^2 + 2Q(x, y)y' + P(x, y)$, {http://dx.doi.org/10.1134/S0012266107050035}, {\textit{Differential Equations}}, \textbf{43}(5), (2007),  595-604.

\bibitem{Bagderina1} Yulia Yu. Bagderina, Invariants of a family of scalar second-order
ordinary differential equations, {doi:10.1088/1751-8113/46/29/295201}, {\textit{J. Phys. A:Math. Theor.}}, \textbf{46}, (2013), 1-36
\end{thebibliography}
\end{document}